\documentclass[a4paper,11pt]{article}
\usepackage{pstricks}
\usepackage{amsmath}
\usepackage{amsfonts}
\usepackage{amscd}
\usepackage{amssymb}
\usepackage[all]{xy}

\newcommand{\al}{{\alpha}}

\newcommand{\om}{{\omega}}

\newcommand{\si}{{\sigma}}

\newcommand{\qed}{\hfill {\bf QED}\medskip}
\newcommand{\Rar}{\longrightarrow}
\newcommand{\ra}{\rightarrow}

\newcommand{\ifff}{if and only if\ }

\newcommand{\MS}{{\medskip}}
\newcommand{\BS}{{\bigskip}}
\newcommand{\NI}{{\noindent}}

\newcommand{\QED}{\hfill {\bf QED}\medskip}

\newtheorem{theorem}{Theorem}[section]
\newtheorem{thm}[theorem]{Theorem}
\newtheorem{cor}[theorem]{Corollary}

\newtheorem{rem}[theorem]{Remark}
\newtheorem{lemma}[theorem]{Lemma}
\newtheorem{prop}[theorem]{Proposition}

\newtheorem{con}[theorem]{Conjecture}
\numberwithin{equation}{section}

\title{\bf Existence of symplectic structures on torus bundles
over surfaces}
\author{Rafa\l\ Walczak \thanks{The author
was partially supported by Grant 2 P03A 036 24 of Polish Committee
of Sci. Research}}

\date{\today}

\begin{document}

\maketitle

\begin{abstract}
Let $E$ be the total space of a locally trivial torus bundle over
the surface $\Sigma_g$ of genus $g>1$. Using the Seiberg--Witten
theory and spectral sequences we prove that $E$ carries a
symplectic structure if and only if the homology class of the
fiber $[T^2]$ is nonzero in $H_2(E,\mathbb R).$

\smallskip


\noindent {\bf AMS classification(2000):} 53D05, 57R57, 55R20

\end{abstract}

\tableofcontents

\section{Introduction}\label{into}

Let $E$ be a closed 4--manifold. The problem whether $E$ admits a
symplectic structure is in general very difficult. There are,
however, some theorems on existence of symplectic forms, for
example Thurston's construction \cite{TH1}. Recalling that
$\Sigma_g$ denotes the surface of genus $g,$ his theorem can be
stated as follows:
\begin{thm}\label{thc1} {\bf [Thurston]} Let $E
{\buildrel \pi \over \longrightarrow} \Sigma_g$ be a surface
$\Sigma_h$--bundle over a surface $\Sigma_g.$ If the homology
class of the fiber is nonzero in $H_2(E,\mathbb R),$ then $E$
admits a $\pi$--compatible (i.e. compatible with the bundle
structure) symplectic form.
\end{thm}

Consider a locally trivial surface bundle over a surface. If the
fiber is different from a torus, Thurs\-ton's construction gives a
symplectic form on $E$. If the fiber is the torus, the situation
is more complicated. If the base is the sphere $S^2$, then all
torus bundles are principal and only the trivial bundle admits a
symplectic structure, as for nontrivial bundles we have $b_2=0$
(by $b_2$ we mean the dimension of the second cohomology group
$\dim H^2(*,\mathbb R)$). If the base is the torus, then Geiges
(\cite{G}) showed that the total space always supports a
symplectic form. The case when the base is a surface of genus
$g>1$ was still unsolved in general.

The principal motivations for us were relations between the space
of all symplectic forms and the subspace of invariant symplectic
forms, assuming $M$ is endowed with a free circle action. In
\cite{HW} we were concerned with the space $S_{inv}$ of such
invariant symplectic forms on 4--manifolds. Here we study torus
fibration from viewpoint of existence of free circle action
preserving fibers. This vehicle recognizes these bundles which
admit symplectic structure. That's why in section \ref{main} we
determine necessary and sufficient condition for existence of such
action (Lemma \ref{id}). We treat the two cases (i.e. these
admitting and not admitting such actions) separately. In the case
where there is no such action we use the {\it Leray--Serre}
spectral sequences in order to prove \ref{gh}. In the other case
we use relevant facts from \cite{HW} and results from section
\ref{tp} to achieve our goals.

Another motivation for us was the following (see \cite{B2}):

\begin{con} [Baldridge] \label{tn} Let $E$ be a closed oriented
4--manifold admitting a symplectic structure and equipped with a
free action of the circle $S^1.$ Then the quotient $E \slash S^1$
fibers over $S^1$. \end{con}

Conjecture \ref{tn} is a generalization of Taubes' Conjecture,
which we quote below.

\begin{con} [Taubes] \label{ta} Let $M$ be a closed oriented 3--manifold.
Then $S^1 \times M$ admits a symplectic structure \ifff $M$ fibers
over $S^1$.\end{con} Note that Taubes' Conjecture \ref{ta} is
equivalent to the existence of invariant symplectic form under the
action of $S^1$ on the first factor of $S^1 \times M^3.$

Total spaces of torus bundles equipped with free circle action
preserving fibers over surfaces form a natural class to check the
conjecture \ref{tn}. We prove it for this class in Section
\ref{tp} by showing that whenever $E$ is a torus bundle over
$\Sigma_g$ with genus greater than 1 (here and in the sequel
$\Sigma_g$ will denote the oriented surface of genus $g$), which
admits a free circle action preserving fibers and a symplectic
form, then $E\slash S^1 \cong \Sigma_g \times S^1.$

\bigskip

The main result of this paper is the following

\bigskip

\NI {\bf Theorem \ref{gh}} {\it The total space $E$ of the
fibration $T^2 \ra E {\buildrel \pi \over \longrightarrow}
\Sigma_g \ \ (g>1)$ is symplectic if and only if the homology
class $[T^2]$ represented by the fiber is nonzero in
$H_2(E,\mathbb R).$}

\bigskip

We should recall that there exists a topological characterization
of 4--dimensional symplectic manifolds (\cite{Go}). They can be
characterized by existence of the {\it hyperpencil}. This
characterization, however, does not seem to be applicable in the
case of our interests, i.e. torus bundles over the surface
$\Sigma_g$ of genus $g>1$. The reason is that there is no evident
proof of nonexistence of the hyperpencil on these torus bundles
which do not admit any symplectic structure.

The main tools used in this paper are the Seiberg--Witten
invariants and spectral sequences.

Recently Etg\"{u} proved in \cite{E} that a principal torus bundle
$E$ over $\Sigma_g$ admits no symplectic form if $g>1$ provided
that $E$ is a product of the circle $S^1$ and a three--dimensional
manifold which is a nontrivial $S^1$ bundle over $\Sigma_g \
(g>1).$ His result verifies positively Taubes' Conjecture \ref{ta}
for manifold $M$ which are $S^1$--principal bundles over $\Sigma_g
\ (g>1).$ In the first part of this paper we prove a theorem which
extends his results. Before stating the theorem let us recall
basic definitions concerning principal torus bundles and the Euler
class of such bundles.

Each $T^2$--principal bundle is obtained by pulling back the
universal bundle $ET^2\rightarrow BT^2$. Since $BT^2$ is homotopy
equivalent to $BS^1 \times BS^1$, therefore $T^2$--principal
bundles over a surface $\Sigma_g$ are classified by a pair $(m,n)$
of two integers. This pair is also called the {\it Euler class} of
the bundle. The Euler class has yet another interpretation. Take
the product $\Sigma_g \times T^2$ and choose any disc $D^2
\hookrightarrow \Sigma_g.$ Remove the counterimage $\pi^{-1}(D^2)$
and glue it back via the identification map $T^2 \times \partial
D^2 \ra T^2 \times \partial D^2$ given by the formula
$((x_1,x_2),\theta) \mapsto ((x_1+\frac{m\theta}{2\pi},
x_2+\frac{n\theta}{2\pi}),\theta).$ The resulting bundle has the
Euler class equal to $(m,n).$

In Section \ref{tp} we consider principal $T^2$--bundles over
surfaces and prove the following:

\bigskip

\NI {\bf Theorem \ref{g}} {\it Let $T^2 \rightarrow E {\buildrel
\pi \over \longrightarrow} \Sigma_g$ be a principal fibration over
the surface $\Sigma_g$, where $g>1$. Let $(m,n)$ be the Euler
class of the fibration. If $mn \neq 0$, then $E$ is not a
symplectic manifold.}

\bigskip

Theorem \ref{g} is proved in two steps. In the first we prove,
following  Etg\"{u} \cite{E}, that $ 0 \in H^2(E,\mathbb Z)$ can
be the only basic class. In the second we write the formula
counting the SW--invariant of $ 0 \in H^2(E,\mathbb Z)$ and
justify that $\underline{sw}_E^4(0)$ is always even, thus the
conclusion comes from Taubes' theorem.

Theorem \ref{g} generalizes directly to the case of torus bundles
which are circle bundles over total spaces of nontrivial circle
bundles over surfaces $\Sigma_g, \ g>1.$ The only difference is
that the formula for $\underline{sw}_E^4(0)$ slightly changes, but
the proof remains the same.

\bigskip

The second part is devoted to arbitrary (locally trivial and
orientable) torus bundles over $\Sigma_g, \ g>1.$ In this part we
prove Theorem \ref{gh}. In fact our proofs give a method to decide
whether a given bundle satisfy the condition of Theorem \ref{gh}
(i.e. when $[T^2] \neq 0$ in $H^2(E,\mathbb Z)$). This theorem
completes the classification of surface bundles over surfaces
admitting symplectic structures. Theorem \ref{gh} is equivalent to
the fact that the total space $E$ of the fibration admits a
symplectic structure \ifff it supports a $\pi$--compatible
symplectic form (\cite{TH1,MS}). Obviously, the only nontrivial
implication is that whenever $E$ supports a symplectic structure,
it has a $\pi$--compatible symplectic form.

\bigskip

We describe now this paper in some details. In section \ref{tp} we
are mainly concerned with principal $T^2$--bundles over $\Sigma_g
\ (g>1).$

\smallskip

We start with an observation that only $0 \in H^2(E,\mathbb Z)$
can be basic class. In order to prove this we use similar
reasoning as in \cite{E}. Next we write the formula for
SW--invariant of $0 \in H^2(E,\mathbb Z)$ using \cite{B1,B2,MOY}.
In order to use Taubes' theorem \ref{t} we prove that the formula
always gives an even number. In Section \ref{int} are gathered
preliminary notions and facts used in the proof of Theorem
\ref{gh}. Further details of the proof are contained in Section
\ref{tp}.

\bigskip

In section \ref{main} we deal with general torus bundles over
surfaces and we finish the proof of the Theorem \ref{gh}. We
consider two cases:

\begin{enumerate}
\item bundles which does not admit any free circle action
preserving fibers, \item bundles which support such action.
\end{enumerate}

A key idea here is that the property that $E$ admits a free circle
action preserving fibers can be described in terms of monodromy
(Proposition \ref{per}). This allows us to prove a useful
trichotomy (Lemma \ref{id}) which gives a connection between the
Betti number $b_1(E)$ and existence of such actions for flat
fibrations. We next prove that if $E$ does not admit such actions,
then we have $[T^2] \neq 0$ and $E$ is symplectic (Theorem
\ref{thc}). If $E$ supports such actions, then we prove that $E$
is symplectic \ifff $E\slash S^1 \cong \Sigma_g \times S^1$ and
$E$ is not a principal torus bundle over $\Sigma_g.$ We also prove
Conjecture \ref{tn} in the case where $E \slash S^1$ is a circle
bundle over $\Sigma_g,\ g>1.$

\bigskip

To prove Theorem \ref{gh} we also apply the Leray--Serre spectral
sequence of the fibration $T^2 \ra E {\buildrel \pi \over
\longrightarrow} \Sigma_g$ to establish two facts. The first
states that $E_{11}^\infty$ depends only on monodromy and not on
the Euler class of the bundle. The second (Lemma \ref{sps}) gives
two conditions $E_{02}^\infty \cong \mathbb Z,$ and $E_{20}^\infty
\cong \mathbb Z$, which both are equivalent to the fact that
$[T^2] \neq 0.$

These facts give answer in Case 1. To see that consider first a
flat bundle $E.$ Corollary \ref{fl} gives that $[T^2] \neq 0,$
which yields that
$$b_2(E)=2+rank(E_{11}^\infty).$$

Consider next the bundle $E'$ which has the same monodromy, but
possibly different Euler from $E.$ Due to Lemma \ref{id} we have
then $b_2(E)=b_2(E')$. Since
$rank(E_{11}^\infty)=rank({E'}_{11}^\infty)$ we have
$$b_2(E')=2+rank({E'}_{11}^\infty).$$

This implies (Lemma \ref{sps}) that $[T^2] \neq 0$ and means that
all bundles in case 1 admit a $\pi$--compatible symplectic form.

\smallskip

Our approach to Case 2 is different. From section \ref{tp} we know
that if $E$ is symplectic, then $E\slash S^1 \cong \Sigma_g \times
S^1.$ Using results form \cite{HW} we prove that if $E$ is not a
principal fibration (i.e. $E$ has nontrivial monodromy), then $E$
admits a symplectic from. This is done in the following manner.

Recall that in \cite{HW} we gave necessary and sufficient
conditions for existence of invariant symplectic forms on
four--dimensional manifolds $E$ equipped in a free circle action.
Following Bouyakoub (\cite{Bo}) we define the subspace

$$L = \{ \alpha \in H^1(E\slash S^1,\mathbb R) \mid
\alpha \cup c_1(\xi) = 0\},$$

\NI where $\xi$ is the principal fibration $S^1 \ra E {\buildrel
\xi \over \longrightarrow} E\slash S^1.$ Then $E$ supports an
invariant symplectic form \ifff $L$ possesses a cohomology class
which has a closed and nondegenerate representative.

For $E\slash S^1 \cong \Sigma_g \times S^1$ we obtain that only
the classes which are pullbacks from $\Sigma_g$ do not have such
representatives. But for these pullback classes we have
$$L=\{\pi_1^*a \mid a \in H^1(\Sigma_g, \mathbb R)\},$$
\NI (where $pr_1: \Sigma_g \times S^1  \ra \Sigma_g$ is the
projection) which means that the bundle $E$ is principal with the
Euler class $(*,0).$ For all other cases we have that $E$ admits
an invariant symplectic structure and $[T^2] \neq 0$ (see
\ref{com}).

In \ref{last} we decide which torus bundles $E$ have $[T^2] = 0$
in terms of monodromy and the Euler class.

Finally we give a couple of applications of Lemma \ref{av} which
concerns averaging of $\pi$--compatible symplectic forms. They are
given in Remarks \ref{avr} and \ref{avrr}.

\bigskip

The author would like to cordially thank Bogus\l\/aw Hajduk, Jan
Dymara and Jarek K\c edra for many hints, suggestions and
explanations. I would like also to thank Bogus\l\/aw Hajduk for
help in editing this paper.

\section{$Spin^c$--structures and SW invariants}\label{int}

In this section we recall some definition and theorems on
$Spin^c$--structures and SW invariants that we will be using. In
this context we also collect facts concerning circle bundles over
3--dimensional manifolds.

\subsection{Classification and pulling back of
$Spin^c$--structures}\label{do}

The set of $Spin^c$--structures on 4--dimensional closed, oriented
manifold $E$ is classified by $H^2(E,\mathbb Z),$ but no such
correspondence can be chosen canonical. Indeed, given two
$Spin^c$--structures $\xi_1,\xi_2$ on $E$ their difference
$\delta(\xi_1,\xi_2)$ is a well--defined element of $H^2(E,\mathbb
Z).$ To see this consider the following diagram.

\centerline{\xymatrix{BS^1 \ar[r] \ar[d]& BS^1 \ar[d] \\
                      E' \ar[r] \ar[d] & BSpin^c(4) \ar[d]\\
                      E \ar@/^1pc/[u]^{\sigma_2} \ar[r]^f \ar[ru]^{f_2} & BSO(4)}}
\smallskip
\NI In this diagram $f$ is the classifying map for $TE$ and $f_2$
is the mapping classifying $Spin^c$--structure $\xi_2.$
Furthermore, bundle $\Xi=(BS^1,E',E)$ over $E$ is the pullback of
$BS^1$--bundle over $BSO(4).$ Since $S^1$ is abelian we have that
$BS^1$ is an H--group and therefore $\Xi$ can be consider as a
principal $BS^1$--bundle over $E.$ Moreover, section $\si_2$ of
$\Xi$ corresponding to $f_2$ gives a trivialization of this
bundle. This trivialization gives a bijection between
$Spin^c$--structures on $E$ and homotopy classes $[E,BS^1]$ of
section of $\Xi,$ which are naturally isomorphic to $H^2(E,\mathbb
Z).$ Once we have fixed $\xi_2$ we define
$$\delta(\xi_1,\xi_2):=\xi_1 - \xi_2 \in H^2(E,\mathbb Z)$$ as element
corresponding to $\xi_1.$ This element is well--defined regardless
of which $Spin^c$--structure we use to trivialize $\Xi.$

\smallskip

Here and subsequently the set of all $Spin^c$--structures on $E$
will be denoted by $\mathcal S(E).$ Recall (\cite{TA3}) we have a
map $$c_1:\mathcal S(E) \ra H^2(X,\mathbb Z)$$ which assigns to a
$Spin^c$--structure $\xi \in \mathcal S(E)$ its first Chern class.
Observe that assignment $$\xi_1 - \xi_2 \mapsto
c_1(\xi_1)-c_1(\xi_2)$$ yields a mapping
$$H^2(E,\mathbb Z) \ni x \mapsto 2x \in H^2(E,\mathbb Z).$$ This
is a direct consequence of the following fact.

\begin{prop}\label{ea} Under the above assumptions
$$c_1(\xi_1)-c_1(\xi_2)=2\delta(\xi_1,\xi_2).$$
\end{prop} For the proof of Proposition \ref{ea} the reader is
referred to \cite{FM}.

\smallskip

Proposition \ref{ea} gives that if $H^2(E,\mathbb Z)$ is $\mathbb
Z_2$--torsion free, then the first Chern class associated to a
$Spin^c$--structure gives an injection of $S(E)$ in $H^2(E,\mathbb
Z).$ However, in presence of $\mathbb Z_2$--torsion there is still
a natural choice of unique element $a$ satisfying
$c_1(\xi_1)-c_1(\xi_2)=2a,$ since we have
$c_1(\xi_1)-c_1(\xi_2)=2\delta(\xi_1,\xi_2).$ This assignment also
classifies $Spin^c$--structures with equal first Chern classes of
their determinants line bundles. They are classified by $A \subset
H^2(E,\mathbb Z),$ where

\begin{eqnarray} A = \{a \in H^2(E,\mathbb Z) \mid 2a=0\}.
\end{eqnarray}

\NI Analogous facts hold in 3--dimensional case. For more details
see \cite{GS}, chapter 10.

\bigskip

Our task is to write an explicit formula for
$\underline{sw}_E^4(0)$ and prove that it always gives even
number. We will to this by applying Baldridge's formula \ref{b2}.
In order to apply this formula we need to establish the set of
these $Spin^c$--structures $\xi$ such that $c_1(\pi^*(\xi))=0.$

We know (\cite{B2}) that if $\pi:E \ra M$ is a principal
$S^1$--bundle and $\xi$ is a $Spin^c$--structure on $M,$ then
$c_1(\pi^*(\xi))=\pi^*(c_1(\xi)).$ The other pulled back
$Spin^c$--structures are obtained by adding of classes $\pi^*(a)$
for $a \in H^2(M,\mathbb Z).$

\smallskip

We will determine the elements $a \in H^2(M,\mathbb Z),$ for which
$\pi^*(a)=0.$ To do this apply the Gysin sequence to get the
following exact sequence $$\cdots \ra H^0(M,\mathbb Z) {\buildrel
\cup\chi \over \longrightarrow} H^2(M,\mathbb Z) {\buildrel \pi^*
\over \longrightarrow} H^2(E,\mathbb Z) \ra H^1(M,\mathbb Z) \ra
\cdots,$$ where $\chi$ denotes the Euler class of the bundle
$\pi:E \ra M.$ It yields that kernel of the mapping $\pi^*:
H^2(M,\mathbb Z) \ra H^2(E,\mathbb Z)$ is equal to the subgroup
generated by the Euler class $\chi.$

This simple fact allows us to determine all $Spin^c$--structures
on $M$ whose pullback has vanishing first Chern class. They can be
described by the condition that their first Chern classes belong
to subgroup $\langle\chi\rangle < H^2(M,\mathbb Z)$ generated by
the Euler class of the fibration.

\bigskip

Next we apply these results to case where $M$ is a nontrivial
principal $S^1$--bundle $\eta'$ over $\Sigma_g, \ (g>1)$ whose
Euler class $e$ equals $n[\Sigma_g]$ for some $n \neq 0.$

It is relatively easy to compute $H^2(M,\mathbb Z).$ We start with
recalling the formula for the group $H_1(M,\mathbb Z).$ We thus
get

\begin{eqnarray}\label{samp} H_1(M,\mathbb Z) \cong \ \
\langle b_1,b_2, \cdots, b_{2g}, x \mid  nx \rangle \cong \mathbb
Z^{2g} \oplus \mathbb Z_n,\end{eqnarray} where $x$ denotes the
homology class of the fiber and $b_1,b_2, \cdots, b_{2g}$ the
generators of the base. Here and subsequently $\langle\cdots \mid
\cdots\rangle$ will denote the abelian group given by generators
and relations.

By Poincare duality we have $PD: H_1(M,\mathbb Z) \cong
H^2(M,\mathbb Z).$ We also have (see \cite{BT}, Ch. 1, Sec. 6)
that $PD(x)=e.$ This yields that $H^2(M,\mathbb Z)$ decomposes as
$T \oplus \mathbb Z^{2g},$ where $T$ denotes finite group
generated by $e.$

This yields that for $n$ odd group $H^2(M,\mathbb Z)$ does not
have 2--torsion and the first Chern class gives an injection of
$Spin^c$--structures in $\langle e \rangle < H^2(M,\mathbb Z).$
For $n$ even there are two elements of order two in $H^2(M,\mathbb
Z),$ therefore there are precisely two $Spin^c$--structures
corresponding to the same first Chern class. Similarly, these
first Chern classes belong to $\langle e \rangle.$

\smallskip

This will determine all $Spin^c$--structure $\xi$ on $M$ such that
$c_1(\pi^*(\xi))=0.$ To do this let us introduce the following
definition.

For $m \in \mathbb Z, n \in \mathbb N$ such that $mn \neq 0$
define

$$A_{m,n}=\{x \in \mathbb Z_n \mid (\exists k \in \mathbb Z) \ x
\equiv km (\hskip-10pt \mod n)\},$$ where $(\mathbb Z_n,+)$
denotes the additive group $\{0,\cdots,n-1\}$ modulo $n.$ Note
that $A_{m,n}$ is a subgroup of $\mathbb Z_n$ generated by element
$m (\hskip-8pt \mod n).$

\smallskip

We have that for $n \in \mathbb Z$ the set of first Chern classes
of these $Spin^c$--structures $\xi$ such that $c_1(\pi^*(\xi))=0$
corresponds to $A_{m,|n|} < T \cong \mathbb Z_n.$

\subsection{Seiberg--Witten invariants}\label{sw}

In this section we review facts and theorems on for the
Seiberg--Witten invariants on three and four--dimensional
manifolds which we will be using.

A powerful tool to deal with the problem of existence of
symplectic structures is the Seiberg--Witten theory. In 1994
Taubes (\cite{TA1},\cite{TA2}) proved the following strong theorem
which yields obstruction to existence of symplectic forms. Namely,
for any symplectic form there exists a {\it compatible} almost
complex structure $J,$ for which we can in turn canonically
associate a $Spin^c$--structure $\xi$ (more details can be found
in \cite{LM},\cite{MS}). This $Spin^c$--structure has a
Seiberg--Witten invariant $SW(\xi) \in \mathbb Z.$

\begin{thm}\label{t} {\bf[Taubes]} Let $X$ be a closed 4--manifold
with $b_+>1$ and a symplectic form $\omega$. Then there is a
canonical $Spin^c$--structure $\xi$ on $X$ such that
$SW_X^4(\xi)=\pm 1$ and $\al=c_1(\det(\xi))$ is the canonical line
bundle $K$ of $(X,\omega)$. Moreover, if
$\underline{sw}_X^4(\alpha) \neq 0$ (i.e. $\alpha$ is a basic
class), then
$$|\alpha \cdot [\omega]|\leq |c_1(K) \cdot [\omega]|$$ and
$\alpha \cdot [\omega]=0$ if and only if $\alpha=0$; $|\alpha
\cdot [\omega]|=|c_1(K) \cdot [\omega]|$ if and only if
$\alpha=\pm c_1(K)$.  \end{thm}

\NI The theorem \ref{t} can be applied to decide whether $E$ is
symplectic in the following way. If all nonzero values of SW
invariants for $Spin^c$--structures are different from $\pm 1$,
then $E$ admits no symplectic structure. These $Spin^c$--
structures (or equivalently elements of $H^2(X,\mathbb Z)$) with
nonvanishing SW--invariants are called {\it basic.} For example,
$\#3\mathbb C \mathbb P^2$ admits no symplectic structure although
it possesses almost complex structure and positive second Betti
number.

\smallskip

Let us now recall the notion of the {\it Seiberg--Witten
polynomial.} As we said in Section \ref{do} there are two ways of
associating $Spin^c$--structures with $H^2(X,\mathbb Z).$ The
first uses the fibration $BS^1 \ra BSpin^c(4) \ra BSO(4)$ induced
from sequence of homomorphisms $S^1 \ra Spin^c(4) \ra SO(4).$ In
this case the set $\mathcal S(X)$ of all $Spin^c$--structures is
in bijective correspondence with $H^2(X,\mathbb Z),$ but no such
correspondence is natural without choosing first a distinguished
element in $\mathcal S(X).$ However, there is the first Chern
class map
$$c_1:\mathcal S(X) \ra H^2(X,\mathbb Z)$$ which assigns
to a $Spin^c$--structure $\xi$ its canonical class $K.$ This class
is induced by the canonical homomorphism $Spin^c(4) \ra U(1).$
This gives another map of $\mathcal S(X)$ to $H^2(X,\mathbb Z).$
This mapping is injective if there is no 2--torsion in
$H^2(X,\mathbb Z)$ (see \cite{LM}, Appendix D).

When $X$ is compact, oriented 4--manifold with $b_2^+>1,$ then the
Seiberg--Witten invariants define, via the map $c_1,$ a map

\begin{eqnarray} \underline{sw}_X^4:H^2(X,\mathbb Z) \ra \mathbb Z
\end{eqnarray}
\NI which is defined up to $\pm 1$ without any additional choices.
That is,

\begin{eqnarray}\label{swm} \underline{sw}_X^4(z) =
\sum_{ \{s \ \mid \ c(s)=z\} } SW(s),
\end{eqnarray}

\NI where $SW(s)$ denotes the value of the Seiberg--Witten
invariant on the class $s \in \mathcal S(X).$

Note that if there is 2--torsion in $H^2(X,\mathbb Z),$ then for
$a \in H^2(X,\mathbb Z)$ we can calculate $\underline{sw}_X^4(a)$
by the same formula \ref{swm}.

It proved useful to package the map $\underline{sw}_X^4$ in a
manner which will now be described. We starts with the group ring
$\mathbb ZH^2(X,\mathbb Z),$ the free $\mathbb Z$ module generated
by the elements of $H^2(X,\mathbb Z).$ Then the invariant
$\underline{sw}_X^4,$ in the $b_2^+ > 1$ case, can be written as
an element in $\mathbb ZH^2(X,\mathbb Z),$ called the {\it
Seiberg--Witten polynomial} of $X:$

\begin{eqnarray}\label{poly} SW_X^4 =
\sum_{z} \underline{sw}_X^4(z)z.
\end{eqnarray} This formula means that for given $z \in
H^2(X,\mathbb Z)$ its SW invariant is given by
$\underline{sw}^4_X(z).$

More details on Seiberg--Witten invariants can be found in
\cite{TA3}.

\bigskip

Our proof of nonexistence of symplectic structures on principal
torus fibrations $E \ra \Sigma_g,$ for $g>1$ and the Euler class
$(m,n)$ satisfying $mn \neq 0$ begins with the observation that we
have a decomposition of $E$ as circle bundle $\eta$ over
three--dimensional manifold $M \cong E\slash S^1,$ which itself is
a nontrivial circle bundle $\eta'$ over $\Sigma_g.$ This helps us
to prove that $0 \in H^2(E,\mathbb Z)$ can be the only basic class
and to write the formula for SW(0). For the first task we follow
Etg\"{u} (\cite{E}). Namely, from the second part of Taubes'
Theorem \ref{t} we know that a basic class cannot be a nonzero
torsion element. Furthermore, by the results cited below we also
get that any basic class must be torsion since it is included in
pullback of group $T \subset H^2(M,\mathbb Z),$ where $T$ is
generated by the first Chern class $c_1(\eta').$

\begin{thm} [Baldridge \cite{B2}] \label{b2} Let $E$ be a closed
oriented 4--manifold with $b_{+}>1$ and a free circle action. Then
the orbit space $M^3$ is a smooth 3--manifold and suppose that
$\chi \in H^2(M, \mathbb Z)$ is the Euler class of the circle
action on $E$. If $\xi$ is  $Spin^c$ structure on $E$ with
$SW_E^4(\xi) \neq 0,$ then $\xi=\pi^*(\xi_0)$ for some $Spin^c$
structure on $M$ and
\begin{eqnarray}\label{for12} SW_E^4(\xi) \ = \ \sum_{\xi'\equiv \xi_0 \ mod \
 \chi} SW_M^3(\xi').\end{eqnarray}\end{thm}

\begin{thm}[Mrowka, Ozsv\'{a}th, Yu \cite{MOY}] \label{om} Let $M$ be
a $S^1$ prin\-ci\-pal bundle over a surface $B$ of genus $g\geq 1$
with Euler class $n\lambda,$ where $\lambda$ is the (positive)
generator of $H^2(B,\mathbb Z)$. If $n \neq 0$, then all basic
classes of $M$ are in $\{k\pi^*(\lambda) \mid 0\leq k \leq
|n|-1\}$, where $\pi$ is the projection. Moreover, we have
\begin{eqnarray}\label{for11} \underline{sw}^3_M(k \cdot\pi^*(\lambda))=
\sum_{s  \equiv  k \ (mod \ n)} \underline{sw}^3_{S^1 \times B} (s
\cdot pr_2^*(\lambda)), \end{eqnarray} where $pr_2$ is the
projection $S^1 \times B\rightarrow B.$
\end{thm}

\NI The Seiberg--Witten polynomial of 3--manifold $M$ is defined
similarly as in the case of 4--manifolds (see the formula
\ref{poly}). For other necessary background on Seiberg--Witten
theory concerning 3--dimen\-sion\-al manifolds consult
\cite{MOY,B2}.

\smallskip

Theorem \ref{om} together with \ref{for11} gives an explicit
formula for the SW--invariant of $M$ (\cite{B1}).
\begin{thm}\label{b1} (Baldridge \cite{B1}) Let $\pi:M \rightarrow
\Sigma_g$ be the total space of a circle bundle over a genus $g$
surface. Assume $c_1(M)=n\lambda$ where $\lambda \in
H^2(\Sigma_g,\mathbb Z)$ is the generator and $n$ is an even
number $n=2l \neq 0$, then the Seiberg--Witten polynomial of $M$
is
$$SW^3_M=\sum_{i=0}^{n-1}\underline{sw}^3_M(t^i)t^i=$$
\begin{eqnarray}\label{for13} = sign(n)
\sum_{i=0}^{|l|-1}\sum_{k=-(2g-2)}^{2g-2} (-1)^{(g-1)+i+k|l|}
\left(\begin{array}{c} 2g-2  \\ (g-1)+i+k|l|
 \end{array}\right) t^{2i},\end{eqnarray} where $t = \exp(\pi^*(\lambda))$
and the binomial coefficient $\left(\begin{array}{c} p
\\ q
 \end{array}\right)=0$ for $q<0$ and $q>p.$ For the formula where $n$ is
odd, replace $l$ by $n$ and $t^{2i}$ by $t^i$.
\end{thm}

The Formula \ref{for13} should be understood that the
Seiberg--Witten invariant of the class $k\lambda$ equals to
coefficient at $t^k$ of the Seiberg--Witten polynomial $SW^3_M.$
Observe also that we use the variable $t = \exp(\pi^*(\lambda))$
in order to keep multiplicative notation of the Seiberg--Witten
polynomial $SW^3_M.$

Formula \ref{for12} allows us to write down an explicit expression
which computes $\underline{sw}_E^4(0)$ (recall that only $0 \in
H^2(E,\mathbb Z)$ can be basic). We prove that under our
assumptions this formula always gives even number. It yields that
$E$ is not symplectic due to Taubes' theorem \ref{t}.

\section{Principal torus bundles} \label{tp}

In this section we consider principal $T^2$--bundles over surfaces
$\Sigma_g$ of genus $g>1$ as well as $S^1$ bundles over
3--manifolds which are $S^1$ bundles over $\Sigma_g \ (g>1).$ Our
main result of this Section is the following.

\begin{thm}\label{g} Let $T^2 \rightarrow E {\buildrel
\pi \over \longrightarrow} \Sigma_g$ be a principal fibration over
the surface $\Sigma_g$, where $g>1$. Let $(m,n)$ be the Euler
class of the fibration. If $mn \neq 0$, then $E$ is not a
symplectic manifold.
\end{thm}

Assume we have a principal torus fibration $E \rightarrow
\Sigma_g$ with the Euler class equal to $(m,n)$ (for a complete
definition of the Euler class see the beginning of section
\ref{main}). We will prove that if $mn \neq 0,$ then $E$ does not
support any symplectic structure. The case where $mn =0$ was
settled by Etg\"{u} in \cite{E}. Namely, he proved (Lemma 3.8 in
\cite{E}) that a nontrivial $S^1$--bundle $M^3 {\buildrel \pi'
\over \longrightarrow} \Sigma_g \ (g>1)$ cannot fiber over $S^1.$

This fact can be proven by another argument using a Thurston
theorem \cite{TH2}. Recall (\ref{samp}) that for a bundle with the
Euler class $e=n[\Sigma_g], \ n \neq 0,$ we have $H_1(M,\mathbb Z)
\cong Z^{2g} \oplus \mathbb Z_n.$ This easily yields that the
homology group $H_2(M,\mathbb Z) \cong H^1(M,\mathbb Z) \cong
\mathbb Z^{2g}$ has a base which consists of embedded tori. To see
this take any set of embedded circles $S_i$ representing $b_1,b_2,
\cdots, b_{2g}$ and consider their counterimages
$(\pi')^{-1}(S_i).$ It follows that the Thurston norm vanishes
identically on $H_2(M,\mathbb Z)$. Now if $M$ fibers over the
circle we have that the fiber must be the sphere or the torus
(\cite{TH2}). Since we have $4\leq 2g = b_1(M) \leq 3$ the claim
follows.

To finish the proof of Taubes' Conjecture he proves using the
SW--theory that a nontrivial circle bundle over $\Sigma_g \times
S^1$ whose Euler class is $n\pi_2^*[\Sigma],$ does not support any
symplectic structure. Obviously this bundle is the principal
$T^2$--bundle over $\Sigma$ with the Euler class equal to $(n,0).$

\bigskip

Assume we have a principal torus fibration $E \rightarrow
\Sigma_g$ with the Euler class equal to $(m,n)$ (for a complete
definition of the Euler class see the beginning of section
\ref{main}). The total space $E$ can be decomposed as a principal
$S^1$ bundle over 3--manifold $M$ which is itself a
$S^1$--principal bundle over $\Sigma_g$. This decomposition goes
as follows. We first take a $S^1$ principal bundle $S^1
\rightarrow M {\buildrel \pi' \over \longrightarrow} \Sigma_g$
over $\Sigma_g$ with the Euler class equal to $n\lambda$, where
$\lambda \in H^2(\Sigma_g,\mathbb Z)$ is the orientation class.
Next we take a $S^1$ principal bundle $S^1 \ra X {\buildrel \pi
\over \longrightarrow} M$ with the Euler class equal to $m
\cdot{\pi'}^*\lambda$. We want to prove that the total space $X$
of the last fibration is diffeomorphic to $E$. This procedure of
decomposing is an example of some more general pattern.

\medskip

Assume we are given a $G$--principal bundle $\xi: G \ra E
{\buildrel \pi \over \longrightarrow} Y$ over some CW--complex $Y$
and a normal subgroup $H \lhd G$ where the inclusion is denoted by
$i:H \hookrightarrow G.$ We can write $\pi = E {\buildrel \xi''
\over \longrightarrow}E\slash H {\buildrel \xi' \over
\longrightarrow} Y$, where $\xi'$ is a $G\slash H$--principal
bundle over $Y$ and $\xi''$ is a $H$--principal bundle over
$E\slash H.$ Our goal is to find relations between classifying
maps for these three bundles $\xi, \xi'$ and $\xi''.$

\smallskip

Consider the following diagram

\smallskip

\centerline{\xymatrix{ Y \ar[r]^f \ar[dr]_d & BG \ar[d]^{Br} \\
                                & B(G\slash H)   \\  }}

\NI where $r:G \ra G\slash H$ is the quotient map and $f$ is the
classifying map for the bundle $\xi.$ We obviously get that $d$ is
the classifying map for the $G\slash H$--principal bundle $\xi'$
which was obtained from $\xi$ by dividing the total space $E$ by
the free right $H$--action.

\NI To proceed consider another diagram.

$$
\begin{CD}\label{G}
             E           @>>>          EG     \\
             @V{\pi}VV                @VVV            \\
             Y            @>{f}>>       BG     \\
\end{CD}
$$

\NI Let us divide the total spaces $E$ and $EG$ by the free right
$H$--action. From the commutative diagram

$$
\begin{CD}\label{H}
             E\slash H   @>{\bar{f}}>>  EG\slash H \cong BH  \\
             @VVV                @V{Bi}VV            \\
             Y            @>{f}>>       BG     \\
\end{CD}
$$

\NI we get that ${\bar{f}}: E\slash H \ra BH$ is the classifying
map for the principal $H$--bundle $\xi''.$

\smallskip

In our case we have $G=T^2$ and $H=S^1.$ We also have that
$d^*:H^2(BS^1) \ra H^2(\Sigma_g)$ is given by $d^*(c_1)=n
[\Sigma_g].$ Furthermore, the classifying map $\bar{f}:E\slash S^1
\ra BS^1$ is given by $\bar{f}^*(c_1)=m\pi^*[\Sigma_g].$ This can
be seen by considering diagram

$$
\begin{CD}\label{CL}
             E\slash S^1   @>{\bar{f}}>>  ET^2\slash S^1 \cong BS^1  \\
             @VVV                @V{Bi}VV            \\
             \Sigma_g            @>{f}>>       BT^2     \\
\end{CD}
$$
\bigskip

In order to use Taubes' Theorem \ref{t} we need to show that
$b^2_+(E)>1$. This follows form the vanishing of the signature
$\sigma(E)$ and the fact that $b_2(E)>3.$

For the proof of the first fact see \cite{K}. The second is then
quite clear. We know that $b_1(E)\geq 2g$ (see Formula \ref{se}).
The Euler characteristics $\chi(E)=0,$ which gives
$b_2(E)=2b_1(E)-2\geq 4g-2\geq 6,$ thus $b_2^+ \geq 3.$

\bigskip

We will check that only $0 \in H^2(E,\mathbb Z)$ can be basic. The
idea here is to modify slightly the argument in \cite{E}. First
notice that a basic class cannot be a nonzero torsion class. This
is a direct consequence of the second part of Taubes' theorem
(\ref{t}). To see this recall that if $\al$ is basic and $\al
\cdot [\om]=0,$ then $\al=0.$

\NI We also have that basic class of $E$ cannot be non--torsion.
To justify this notice that all nonzero Seiberg--Witten invariants
$SW_M^3$ are included in a group $T$ which was defined as

\begin{eqnarray}\label{s} T=\{k \cdot {\pi'}^*\lambda \mid k \in \mathbb
Z\}.
\end{eqnarray}
This fact follows directly from the Theorem \ref{om}. Due to
Theorem \ref{b2} we get that basic classes of $E$ are included in
the set $\pi^*T.$ This gives our claim since $T$ is a finite
group.

\smallskip

Theorem \ref{b1} also gives an explicit formula to count the
Seiberg--Witten invariant of the zero class. Seiberg--Witten
invariants of $M$ are given as coefficients of the Seiberg--Witten
polynomial exhibited in a very convenient form in Theorem
\ref{b1}. Take any $Spin^c$--structure $\xi_0$ on $M$ such that
$c_1(\pi^*(\xi_0))=0.$ We know (see Section \ref{do}) that these
$Spin^c$--structures are characterized by the condition that
$c_1(\xi_0) \in A_{m,|n|}<T.$ For each such structure $\xi_0$
consider the set
$$\{\xi \mid \xi \equiv \xi_0 \ (\hskip-10pt \mod \chi)\}$$
corresponding to $\{c_1(\xi_0)\}+A_{2m,|n|} \subset T$ via the
first Chern class. Formula for $\underline{sw}_E^4(0)$ has now the
form

\begin{eqnarray}\label{for23} \underline{sw}_E^4(0)=
\sum_{\{\xi \mid  \  2\chi \mid (c_1(\xi)-c_1(\xi_0)) \ \}}
\sum_{\{\xi_0 \mid \  \chi \mid c_1(\xi_0) \ \}} SW_M^3(\xi).
\end{eqnarray}

This means that

\begin{eqnarray}\label{for25} \underline{sw}_E^4(0)=
\sum_{\{j \in \mathbb Z_n \mid \ (i-j) \in A_{2m,|n|}  \}} \left(
\sum_{i \in A_{m,|n|}} \underline{sw}_M^3(t^j) \right).
\end{eqnarray}

If we denote
$$A=\{i \mid 0\leq i \leq |n|-1,\ \ \gcd(|m|,|n|)\mid i\}$$ and
$$A'=\{i \mid 0\leq i \leq |n|-1,\ \ \gcd(|m|,|n|)\mid i, \ 2 \mid
i\},$$

\NI the formula \ref{for25} can be explicitly given by

$$\underline{sw}_E^4(0)=$$ $$=sign(n)\frac{|n|}{\gcd(|m|,|n|)} \cdot $$
\begin{eqnarray}\label{for36} \cdot \left( \sum_{i \in
A}\sum_{k=-(2g-2)}^{2g-2} (-1)^{(g-1)+i+k|n|}
\left(\begin{array}{c} 2g-2  \\ (g-1)+i+k|n|
 \end{array}\right) \right), \end{eqnarray} if $n$ is odd, and

$$\underline{sw}_E^4(0)=$$ $$=sign(n)\frac{|n|}{x \gcd(|m|,|n|)} \cdot$$
\begin{eqnarray}\label{for35} \cdot \left(
\sum_{i \in A'}\sum_{k=-(2g-2)}^{2g-2} (-1)^{(g-1)+i+k|\frac n2|}
\left(\begin{array}{c} 2g-2  \\ (g-1)+i+k|\frac n2|
\end{array}\right) \right),  \end{eqnarray} if $n$ is even, where $x$ equals 1 or 2
depending on whether $\gcd(|2m|,|n|)=\gcd(|m|,|n|)$ or
$\gcd(|2m|,|n|)=2 \gcd(|m|,|n|).$

\smallskip

For $n$ odd this can be proven in the following way. Since
$\gcd(2m,|n|)=\gcd(m,|n|)$ we have
$\{c_1(\xi_0)\}+A_{2m,|n|}=A_{m,|n|} \subset T.$ Thus for fixed $i
\in A_{m,|n|}$ the sums in \ref{for25} are all equal and they are
counted $\frac{|n|}{\gcd(|m|,|n|)}$ times. Furthermore, for fixed
$i$ these sums are obviously equal to
$$sign(n) \sum_{i \in
A}\sum_{k=-(2g-2)}^{2g-2} (-1)^{(g-1)+i+k|n|}
\left(\begin{array}{c} 2g-2  \\ (g-1)+i+k|n|
 \end{array}\right).$$

For $n$ even the reasoning is similar. Observe also that if
$\gcd(|2m|,|n|)=2 \gcd(|m|,|n|),$ then $\frac{|n|}{2
\gcd(|m|,|n|)}$ is an integer number.

\bigskip

Therefore it is enough to examine parity of the following sums:

\begin{eqnarray}\label{for16}\sum_{i \in
A}\sum_{k=-(2g-2)}^{2g-2} (-1)^{(g-1)+i+k|n|}
\left(\begin{array}{c} 2g-2  \\ (g-1)+i+k|n|
 \end{array}\right), \end{eqnarray} if $n$ is odd, and

\begin{eqnarray}\label{for15}
\sum_{i \in A'}\sum_{k=-(2g-2)}^{2g-2} (-1)^{(g-1)+i+k|\frac n2|}
\left(\begin{array}{c} 2g-2  \\ (g-1)+i+k|\frac n2|
\end{array}\right),  \end{eqnarray} if $n$ is even.

The formulas \ref{for15} and \ref{for16} seemed to be complicated.
Using the Mathematica programm I checked that for
$m,n=-20,\ldots,-1,1,\ldots,20$ and $g=2,\ldots 20$ (which is
approximately 30000 cases) the formulas give even numbers. This
convinced me that the sums should be even.

\smallskip

We will show that the sums are always even. Taubes' theorem will
tell us then that $E$ carries no symplectic structure.

In order to continue, we will consider two cases.

\subsection{$n$ odd}

We analyze the parity of the sum \ref{for16} and \ref{for15},
therefore we omit the powers of $(-1)$ in these formulas.

Summands of the sum $$\sum_{i \in
A}\sum_{k=-(2g-2)}^{2g-2} \left(\begin{array}{c} 2g-2  \\
(g-1)+i+k|n|
 \end{array}\right)$$ are in bijective correspondence with the set
$$B=\{ \left(\begin{array}{c} 2g-2  \\ r
 \end{array}\right) \mid r \in C \},$$ where
$$C=\{(g-1)-(2g-2)|n|, \ldots, (g-1)-(2g-2)|n| +u
\gcd(|m|,|n|),\ldots,$$ $$(g-1)+(2g-1)|n| -\gcd(|m|,|n|)\}$$ and
$u$ is a natural number. To see this just change the order of
summing to get
$$\sum_{k=-(2g-2)}^{2g-2}\sum_{i \in
A} \left(\begin{array}{c} 2g-2  \\
(g-1)+i+k|n|
 \end{array}\right)$$

\NI and expand this sum. Thus for given $k \in
\{-(2g-2),\ldots,2g-2\}$ and $i \in A$ we obtain elements
$(g-1)+k|n|,\ldots,(g-1)+(k+1)|n|-\gcd(|m|,|n|)$ of the set $C.$
Since
$$(g-1)-(2g-2)|n|<0,$$ $$(g-1)+(2g-1)|n| -\gcd(|m|,|n|)>2g-2$$ and
$$\left(\begin{array}{c} 2g-2  \\ g-1
 \end{array}\right) \in B,$$ all nonzero members of
$B$ may be represented as a sum of pairs

\noindent $\left\{ \left(\begin{array}{c} 2g-2  \\ r
\end{array}\right) , \left(\begin{array}{c} 2g-2  \\ 2g-2-r
\end{array}\right) \mid r < g-1,\ r \in C  \right\}$ and a single element

\noindent $\left(\begin{array}{c} 2g-2  \\ g-1
\end{array}\right)$. It suffices to notice that
$\left(\begin{array}{c} 2g-2  \\ g-1
\end{array}\right)$ is an even number
provided that $g>1,$ since
$$\left(\begin{array}{c} 2g-2  \\ g-1 \end{array}\right)=2
\left(\begin{array}{c} 2g-3  \\ g-1 \end{array}\right).$$

\subsection{$n$ even}

The case when $n$ is even and $m$ is odd is proved in exactly the
same way as before and will be omitted. Therefore without loss of
generality we can assume that $m$ is also even.

In this case we have $2 \mid \gcd(|m|,|n|),$ so $A=A'$ and
therefore we shall examine the parity of the sum $$\sum_{i \in
A}\sum_{k=-(2g-2)}^{2g-2} \left(\begin{array}{c} 2g-2  \\
(g-1)+i+k|\frac n2| \end{array}\right) ,$$ where $A$ and $A'$ are
defined as in the previous case. Let us divide the sum into two
summands: $$\sum_{i \in A}\sum_{k \in
F} \left(\begin{array}{c} 2g-2  \\
(g-1)+i+k|\frac n2|
 \end{array}\right)$$ and $$\sum_{i \in
A}\sum_{k \in F'} \left(\begin{array}{c} 2g-2  \\
(g-1)+i+k|\frac n2|
 \end{array}\right),$$ where $F=\{-(2g-2),-(2g-4),\ldots,2g-4,2g-2 \}$ and
$F'=\{-(2g-3),-(2g-5),\ldots,2g-5,2g-3 \}.$ Thus the first sum may
be expressed as $$\sum_{i \in A}\sum_{k=-(g-1)}^{g-1} \left(\begin{array}{c} 2g-2  \\
(g-1)+i+k|n|
 \end{array}\right).$$ By the same reasoning as in
the case where $n$ is odd we deduce that this sum is even (observe
that $(g-1)+(-g-1)|n|<0$ and $(g-1)+(g+1)|n|>0$ since $n$ is a
nonzero even number). As for the second sum, similarly as in the
previous case, we notice that elements of this sum are in
bijective correspondence with the set
$$B'=\{ \left(\begin{array}{c} 2g-2  \\
r
 \end{array}\right)  \mid r \in C' \},$$ where
$$C'=\{(g-1)-(2g-3)| \frac n2|, \ldots, (g-1)-(2g-3)|\frac n2| +u
\gcd(|m|,|n|),\ldots,$$ $$(g-1)+(2g-1)|\frac n2|
-\gcd(|m|,|n|)\}$$ and $u$ is a natural number.  Observe that
$$(g-1)-(2g-3)| \frac n2|\leq 0,$$  $$(g-1)+(2g-1)|\frac n2|
-\gcd(|m|,|n|)\geq 2g-2,$$ and that all nonzero members of $B'$
are distributed symmetrically with respect to $\left(\begin{array}{c} 2g-2  \\
g-1
 \end{array}\right)$. To be more precise: if $$y=(g-1)+k |\frac n2| +x
\gcd(|m|,|n|) \in C',$$ where $0\leq x
<\frac{|n|}{\gcd(|m|,|n|)}$, $k \in F'$  and $0 \leq y\leq g-1$,
then

\NI $k \not\in \{-(2g-3),2g-3\},$
$$z=(g-1)-(k+2)|\frac n2| +
(\frac{|n|}{\gcd(|m|,|n|)}-x)\gcd(|m|,|n|)  \in C',$$ $k \in F'$
and $z\geq g-1$ as well, since we have $z+y=2g-2$. \qed

\BS

The theorem we have just proved can be slightly generalized.
Namely, we do not need the Euler class to be included in $T$ (see
\ref{s}). Assume then that $\chi \not\in T.$ Observe that the same
argument as before yields that only zero class can be basic. The
difference is that the formula for $\underline{sw}_E^4(0)$
slightly changes. To be more precise, we have that
$\langle\chi\rangle \cap T=\{0\}$ and therefore $\{\xi\mid \
2\chi\mid(c_1(\xi_0)-c_1(\xi))\}=\{\xi_0\}$ for every
$Spin^c$--structure $\xi_0$ such that $\pi^*(c_1(\xi_0))=0.$ Thus
Formulas \ref{for35} and \ref{for36} have the following forms.

$$\underline{sw}_E^4(0)=$$
\begin{eqnarray}\label{for46}=sign(n) \sum_{i \in
A}\sum_{k=-(2g-2)}^{2g-2} (-1)^{(g-1)+i+k|n|}
\left(\begin{array}{c} 2g-2  \\ (g-1)+i+k|n|
 \end{array}\right), \end{eqnarray} if $n$ is odd, and

$$\underline{sw}_E^4(0)=$$
\begin{eqnarray}\label{for45} =sign(n)
\sum_{i \in A'}\sum_{k=-(2g-2)}^{2g-2} (-1)^{(g-1)+i+k|\frac n2|}
\left(\begin{array}{c} 2g-2  \\ (g-1)+i+k|\frac n2|
\end{array}\right),  \end{eqnarray} if $n$ is even.

These are precisely Formulas \ref{for15} and \ref{for16}. We
already proved they represent even numbers provided that $g>1$ and
$n \neq 0$. \BS

\begin{cor}\label{pr} If $T^2 \ra E \ra \Sigma_g$ is a nontrivial principal
fibration, then homology class of the fiber $[T^2]=0$ in
$H_2(E,\mathbb R).$
\end{cor}

\begin{rem}{\em
Note that in this section we proved our Conjecture \ref{tn} for
torus bundles equipped with free circle action preserving fibers.
Recall that the Conjecture \ref{tn} states that if  $E$ is a
closed 4--manifold equipped with a free circle action admitting a
symplectic structure, then the quotient $E \slash S^1$ fibers over
$S^1$. Indeed, in our case we have that $E \slash S^1$ is a
principal circle bundle $\xi$ over surface $\Sigma_g.$ If this
surface is $S^2$ or $T^2,$ the Conjecture obviously hold. if the
genus $g$ is greater than 1, then $\xi$ must be trivial (recall
that if $n \neq 0,$ then $E$ is not symplectic). In the next
section we settle the question which $S^1$ bundles over $\Sigma_g
\times S^1$ admit symplectic structures.} \hfill $\Box{}$
\end{rem}

\section{General torus bundles} \label{main}

In this section we discuss the general case of $T^2$ bundles over
surfaces $\Sigma_g$ of genus $g>1.$

\smallskip

Each such bundle is, up to isomorphism, uniquely determined by the
monodromy homomorphism $\pi_1\Sigma_g \ra \pi_0Diff^+T^2 \cong
SL(2,\mathbb Z)$, which reveals the structure of the bundle over
1--skeleton of the base, and the Euler class. We can assume that
monodromy of such a bundle are linear automorphisms of $T^2,$ thus
they are given by $A_i \in SL(2,\mathbb Z).$ The monodromy can be
also defined by use of {\it classifying spaces.} Namely, if
$f:\Sigma_g \ra BDiff^+T^2$ is the classifying map for the bundle
and $B\pi :BDiff^+T^2 \ra BSL(2,\mathbb Z)$ is the map associated
to the natural map $\pi:Diff^+T^2 \ra \pi_0Diff^+T^2 \cong
SL(2,\mathbb Z)$, then we can define the mapping $\rho$ by the
formula $\rho=(\pi_1)_*(B\pi \circ f):\pi_1(\Sigma_g) \ra
\pi_1(BSL(2,\mathbb Z)) =SL(2,\mathbb Z).$ This mapping is a
representation $\pi_1(\Sigma_g) \ra SL(2,\mathbb Z)$ which also
determines monodromy associated with $E.$

The Euler class, as an obstruction for a cross section of the
bundle, is a pair of integer numbers $(m,n),$ which may be defined
as follows. Let $f:\Sigma_g \ra BDiff^+T^2$ be the classifying map
for the bundle, and $B(ev) :BDiff^+T^2 \ra BT^2$ be the H--map
associated to the H--homomorphism $ev:Diff^+T^2 \ra T^2$ given by
the evaluation. Since $BT^2 \cong \mathbb C \mathbb P^\infty
\times \mathbb C \mathbb P^\infty,$ the map $B(ev) \circ f$
defines the pair of natural numbers $(m,n)$ which are the Euler
class of the bundle.

We call the bundle {\it flat} if its structural group can be
reduced to a discrete group. In case of torus bundle it means that
the Euler class of the bundle is equal to $(0,0)$ and the
structural group reduces to $SL(2,\mathbb Z).$

A realization of a given Euler class $(m,n)$ for a flat bundle can
be described in the following way. Take a flat bundle and choose
any disc $D^2 \hookrightarrow \Sigma_g.$ Remove the counterimage
$\pi^{-1}(D^2)$ and glue it back via the identification map $T^2
\times \partial D^2 \ra T^2 \times
\partial D^2$ given by the formula $((x_1,x_2),\theta) \mapsto
((x_1+\frac{m\theta}{2\pi}, x_2+\frac{n\theta}{2\pi}),\theta).$

\smallskip

We begin with the following lemma.

\begin{lemma}\label{comi} Any flat $T^2$ bundle over $\Sigma_g$
supports a $\pi$--compatible symplectic form. \end{lemma}

{\bf Proof.} The total space $E$ of the bundle can be described as
the quotient by the following action of $\pi_1(\Sigma_g)$ on
$\mathbb H^2 \times T^2.$ First take representation
$\rho:\pi_1(\Sigma_g) \ra SL_2\mathbb Z$ induced by monodromy
associated with $E.$ Define next

\begin{eqnarray}\label{for2} g(x,y)=(xg^{-1},\rho(g)y),\end{eqnarray}

\NI where on the first coordinate we have the natural action of
$\pi_1(\Sigma_g)$ on the universal covering of $\Sigma_g.$ To
finish notice that the product symplectic form on $\mathbb H^2
\times T^2$ is invariant with respect to \ref{for2}. \QED

An equivalent formulation of the Lemma \ref{comi} is given by the
following corollary.

\begin{cor}\label{fl} For any flat $T^2$--bundle over
$\Sigma_g$ we have that the homology class of the fiber $[T^2]$ is
nonzero in $H_2(E,\mathbb R).$
\end{cor}

\NI The equivalence between Lemma \ref{comi} and Corollary
\ref{fl} is due to Thurston \cite{TH1}:

\begin{thm}\label{thc} Let $E {\buildrel \pi \over \longrightarrow}
\Sigma_1$ be a surface $\Sigma_2$--bundle over surface $\Sigma_1.$
If the homology class of the fiber is nonzero in $H_2(E,\mathbb
R),$ then $E$ admits a $\pi$--compatible symplectic form.
\end{thm}

Recall now the following facts which we shall need in the sequel.

If $T^2 \ra E {\buildrel \pi \over \longrightarrow} \Sigma_g$ is a
flat bundle, then

$$H_1(E) \cong$$

\begin{eqnarray}\label{fi1}\langle b_1,b_2, \cdots, b_{2g}, x_1,x_2 \mid
A_ix_1-x_1,A_ix_2-x_2, i=1,2,\cdots,2g \rangle.\end{eqnarray}

\NI The generators $x_1,x_2$ give the basis for $H_1(T^2,\mathbb
Z)$ and $b_1,b_2, \cdots,b_{2g}$ are the standard basis for
$H_1(\Sigma_g, \mathbb Z).$

If the  bundle $T^2 \ra E {\buildrel \pi \over \longrightarrow}
\Sigma_g$ is not flat, we have

$$H_1(E) \cong \ \ \langle b_1,b_2, \cdots, b_{2g}, x_1,x_2
\mid$$ \begin{eqnarray}\label{se} | mx_1 + nx_2,
A_ix_1-x_1,A_ix_2-x_2, i=1,2,\cdots, 2g, \rangle, \end{eqnarray}

\NI where the pair $(m,n) \in \mathbb Z^{\oplus 2}$ denotes the
Euler class. Note that the pair $(m,n)$ is determined by the basis
$x_1,x_2$ of $H_1(F).$

\bigskip

To proceed we make an observation that we will use occurs here is
that the question whether $E$ admits preserving fibers action can
be answer in terms of the monodromy. Namely, existence of such
action is equivalent to the existence of $x \in H^1(T^2,\mathbb
Z)$ preserved by the monodromy. This is the content of the next
proposition.

\begin{prop}\label{per} The total space $E$ of a fibration $T^2 \ra E
{\buildrel \pi \over \longrightarrow} \Sigma_g$ supports a free
circle action preserving fibers \ifff there exists a nonzero
element $z \in \mathbb Z^2$ such that $A_iz=z$ for
$i=1,\ldots,2g.$
\end{prop}

{\bf Proof.} Assume that there exists a common eigenvector $z \in
\mathbb Z^2$ with eigenvalue 1 for all monodromy. We can assume
that the bundle $E$ is isomorphic with the bundle whose all
monodromy $A_i \in SL(2,\mathbb Z).$ We can also assume that $z$
is primitive, i.e. $\gcd(a,b)=1,$ where $z=(a,b).$ If $\{z,u\}$ is
a basis of $H_1(T^2,\mathbb Z),$ than all monodromy have the
form $\left(\begin{array}{cc} 1 & * \\
0 & 1 \end{array}\right).$ We have then a free circle action
preserving fibers on $E$ given by the subgroup corresponding to
$z.$

\smallskip

Assume now that $E$ admits such circle action. Any effective
circle action on $T^2$ is determined, up to isomorphism, by the
homology class of the orbit. As a result we can choose a basis
$\{z,u\} \in H_1(T^2,\mathbb Z)$ for which all monodromy matrices
$A_i$ have
the form $\left(\begin{array}{cc} 1 & * \\
0 & 1 \end{array}\right).$ \QED

The Proposition \ref{per} allows us to prove the following
trichotomy.

\begin{lemma}\label{id} Let $T^2 \ra E {\buildrel \pi
\over \longrightarrow} \Sigma_g$ be a flat fibration. Then

\begin{enumerate}
\item $E$ is trivial $\Leftrightarrow b_1(E)=2g+2,$ \item $E$ is
nontrivial and admits a free circle action preserving fibers
$\Leftrightarrow b_1(E)=2g+1,$ \item $E$ does not admit any free
circle action preserving fibers $\Leftrightarrow b_1(E)=2g.$
\end{enumerate}
\end{lemma}

{\bf Proof.} The first equivalence follows from cohomological
computation \ref{fi1} and is obvious.

\smallskip

Assume that $b_1(E)=2g+1.$ We will prove that $E$ supports a free
circle action preserving fibers. Our assumption yields that the
group $S<\mathbb Z^2$ defined as
\begin{eqnarray}\label{des} S=\langle A_ix_1-x_1, A_ix_2-x_2,
i=1,2,\cdots,2g \rangle, \end{eqnarray} where
$\langle\ldots\rangle$ denotes the smallest subgroup generated by
given relations, is isomorphic to $\mathbb Z.$ To proceed take any
basis $\{z,u\} \in \mathbb Z^2$ such that $S < \mathbb Z z.$ Thus
we obtain $A_iz-z=k_i z$ for some integers $k_i \in \mathbb Z.$ In
the basis $\{z,u\}$ the matrices $A_i$ have the form
$\left(\begin{array}{cc} k_i+1 & * \\
0 & * \end{array}\right).$ Additionally we cannot have $k_i+1=-1$
since $A_i$ would have the form $\left(\begin{array}{cc} -1 & * \\
0 & -1 \end{array}\right)$ and $b_1(E)=2g$. Finally we have
that $A_i$ are given by $\left(\begin{array}{cc} 1 & * \\
0 & 1 \end{array}\right).$ But this and Proposition \ref{per}
means that $E$ admits a free circle action preserving fibers.

\smallskip

Assume now that $E$ is nontrivial and admits such circle action.
Due to Proposition \ref{per} we have that in some base $\{z,u\}
\in \mathbb Z^2$ the monodromy matrices have the form
$\left(\begin{array}{cc} 1 & * \\
0 & 1 \end{array}\right).$ This yields that $b_1(E)=2g+1.$

\smallskip

The third equivalence follows from the first two equivalences.
\QED

\bigskip

At this point we would like to analyze the Leray--Serre spectral
sequence for the fibration $T^2 \ra E {\buildrel \pi \over
\longrightarrow} \Sigma_g.$ We need this to establish two facts
concerning our bundles. First is to justify that $E_{11}^\infty$
does not depend on the Euler class, but only on monodromy. Second
is to prove a lemma which states that the homology class $[T^2]$
is nonzero in $H_2(E,\mathbb R)$ \ifff any of the conditions:
$E_{02}^\infty \cong \mathbb Z,$ or $E_{20}^\infty \cong \mathbb
Z,$ hold.

To prove the first fact we follow Geiges (\cite{G}). The
$E^2$--page of the Leray--Serre spectral sequence is given by
$$E_{pq}^2 = H_p(\Sigma_g,\mathcal H_q(T^2)),$$ where $\mathcal H$
denotes the system of local coefficients, and the spectral
sequence converges to $H_*(E,\mathbb Z).$

For our purposes it is enough to consider
$$E_{11}^\infty = E_{11}^2 = H_1(\Sigma_g,\mathcal H_1(T^2)).$$

Using the standard action of $\pi_1(\Sigma_g)$ on the universal
cover we get the following free resolution of $\mathbb Z$ over
$\mathbb Z\Pi,$ where $\Pi \cong \pi_1(\Sigma_g).$

\begin{eqnarray} 0 \ra \mathbb Z\Pi \ra \bigoplus_{i=1}^{2g} \mathbb
Z\Pi \ra \mathbb Z\Pi \ra \mathbb Z \ra 0.\end{eqnarray}

Recall that if $\mathcal{R}$ is a projective resolution of
$\mathbb Z$ over $\mathbb Z\Pi$ and $V$ is a $Q$--module, then
$$H_*(Q,V)=H_*(\mathcal{R} \otimes_{\mathbb Z\Pi}V)$$
by definition.

Tensoring the above resolution with $\mathbb Z \oplus \mathbb Z$
over $\mathbb Z\Pi$ gives

\begin{eqnarray} \label{sqna} 0 \Rar \mathbb Z^{\oplus 2} {\buildrel i_*
\over \longrightarrow} \mathbb Z^{\oplus 4g} {\buildrel p_* \over
\longrightarrow} \mathbb Z^{\oplus 2} \ra \mathbb Z^{\oplus 2} \ra
0.\end{eqnarray}

Then $E_{11}^2=ker(p_*) \slash im(i_*).$

This implies that $E_{11}^2=E_{11}^\infty$ depends only on
monodromy and not on the Euler class. Indeed, if we keep the
monodromy and change the Euler class then the sequence \ref{sqna}
remains the same.

\smallskip

Next we establish the second fact, which will shall be using.

\begin{lemma}\label{sps} Let $T^2 \ra E {\buildrel \pi \over
\longrightarrow} \Sigma_g$ be an orientable fibration where $g>1.$
Then the condition $[T^2] \neq 0$ is equivalent to each of the
following two conditions:  $E_{02}^\infty \cong \mathbb Z,$ or
$E_{20}^\infty \cong \mathbb Z.$
\end{lemma}

{\bf Proof.} In the $E^2$--page of the spectral sequence for the
fibration we have $H_0(\Sigma_g,H_2(F))=E_{02}^2 \cong \mathbb Z$
and $H_2(\Sigma_g,H_0(F))=E_{20}^2 \cong \mathbb Z$ since the
local coefficients systems for $H_0(T^2,\mathbb Z)$ and
$H_2(T^2,\mathbb Z)$ are both trivial. It follows that the summand
$E_{02}^\infty$ has the following interpretation: $[T^2] \neq 0
\Leftrightarrow E_{02}^\infty \cong \mathbb Z.$ This can be easily
seen by noticing that (see for example \cite{Sp}, ch.9, sec.3, p.
482) the homomorphism $i_*:H_2(T^2,\mathbb Z) \ra H_2(E, \mathbb
Z)$ is the composition

$$H_2(T^2,\mathbb Z) \cong H_0(\Sigma_g,H_2(T^2,\mathbb Z)) \cong
E_{02}^2 \ra E_{02}^\infty =$$ \begin{eqnarray} = F_0H_2(E,\mathbb
Z) \subset H_2(E,\mathbb Z),
\end{eqnarray}

\NI where $F_i$ is the filtration (see also \cite{S} and
\cite{TH1}).

Similarly the summand $E_{20}^\infty$ has an analogous
interpretation: $[T^2] \neq 0 \Leftrightarrow E_{02}^\infty \cong
\mathbb Z.$ This can proven the same way as before. The local
coefficients systems for $H_0(T^2,\mathbb Z)$ is trivial, so
(\cite{Sp}, ch.9, sec.3, p. 483) $\pi_*$ is the composition

$$H_2(E,\mathbb Z) = F_2(H_2(E,\mathbb Z)) \ra E_{20}^\infty \ra
E_{20}^2 \cong$$ \begin{eqnarray}\cong H_2(\Sigma_g,
H_0(T^2,\mathbb Z)) \ra H_2(\Sigma_g,\mathbb Z)
\end{eqnarray}

The isomorphism $E_{02}^\infty \cong \mathbb Z$ is equivalent to
the fact that $\pi_*:H_2(E) \ra H_2(\Sigma_g)$ is rationally onto.
This in turn means that $\pi^*[\Sigma_g] \neq 0,$ where
$[\Sigma_g]$ is the orientation class. Since
$PD(\pi^*[\Sigma_g])=[T^2]$ in cohomology with real coefficients
the assertion is proven. \QED

\bigskip

To continue we review some relevant facts from
\cite{Bo,HW,MD,MT,TH1,TH2}.

We denote by $M$ a compact oriented smooth manifold with a smooth
free action of $S^1$ and $N \cong M \slash S^1.$ The space of
invariant symplectic forms consistent with the given orientation
will be denoted by $S_{inv}.$

Let $\pi:M \ra N$ denote the principal $S^1$ -- fibration given by
the action. By $\al$ we will denote the non\-dege\-ne\-rate
(=no\-where va\-ni\-shing) closed 1--form $\alpha$ on $N$
satisfying
$$\pi^*\alpha = \iota_X\omega, $$
where $X$ is the infinitesimal generator of the action.

\begin{lemma}\label{bu} {\bf \cite{Bo}} If $\omega \in S_{inv}$, then
\begin{eqnarray}\label{fro} [\alpha] \cup c_1(\pi) = 0.\end{eqnarray}
\end{lemma}

{\bf Proof}. Take any connection form $\eta \in \Omega^1(M,\mathbb
R)$. By Chern -- Weyl there exists a closed 2 -- form $c_1 \in
H^2(N,\mathbb R)$ such that $\pi^*c_1 = d\eta$ and $c_1$
represents $c_1(\pi) \in H^2(N^3,\mathbb R)$. There also exists a
unique 2 -- form $\beta' \in \Omega^2(N^3,\mathbb R)$ such that
$\omega - \eta \wedge \iota_X\omega = \pi^*\beta'$. By
differentiating both sides of the last equation we get $d\eta
\wedge \iota_X\omega = \pi^*d\beta'$ and this implies the lemma.
\QED

\MS

From now on we assume that $\dim M=4.$

\MS

Define the subspace $L \subset H^1(N^3,\mathbb R)$ by

\begin{eqnarray}\label{na} L = \{ \alpha \in H^1(N^3,\mathbb R) \mid
\alpha \cup c_1(\pi) = 0\}. \end{eqnarray}

Let us also recall the notion of {\bf Thurston norm} \cite{TH2}
(see also \cite{MT}) for a 3--manifold. If $N^3$ is a compact,
connected and oriented manifold without boundary then for any
compact oriented n--component surface $S = S_1 \sqcup \cdots
\sqcup S_n$ embedded in $N$ define

\begin{eqnarray} \chi_-(S) = \sum_{\chi(S_i) <0} \arrowvert
\chi(S_i) \arrowvert. \end{eqnarray}

The {\bf Thurston norm} on $H_2(N^3,\mathbb Z)$ and, by the
Poincare duality, on $H^1(N^3,\mathbb Z)$ is given by

\begin{eqnarray} {\Arrowvert \phi \Arrowvert}_T = \inf \{ \chi_-(S)
\mid [S] = \phi  \}.\end{eqnarray}

The Thurston norm can be extended linearly to $H^1(N^3,\mathbb
R)$. Let $B_T = \{\phi \ : \  {\Arrowvert \phi \Arrowvert}_T \leq
1 \}$ denote the unit ball in the Thurston norm. It is a (possibly
noncompact) polyhedron in $H^1(N^3,\mathbb R)$. Suppose $\phi' \in
H^1(N^3,\mathbb Z)$ is represented by a fibration $N^3 \rightarrow
S^1$. Then $\phi'$ is contained in the open cone $\mathbb R_+
\cdot F$ over a top -- dimen\-sional face $F$ of the Thurston norm
ball $B_T$. In this case we say $F$ is a {\it fibered face} of the
Thurston norm ball. The Thurston norm can be also defined if
boundary of $N^3$ is a union of tori.

Assume now that $N^3$ fibers over the circle, $b_1=m$ and it is
not $S^2 \times S^1$. In this case we have that $L \cap \mathbb
R_+ \cdot F$ has a homotopy type of a point, or a sphere $S^{m-1}$
if $c_1=0$ in $H^2(N,\mathbb R),$ and $S^{m-2}$ when $c_1 \neq 0$
provided that Thurston norm vanishes identically. If $N \cong S^2
\times S^1$, then total space of any circle fibration over $N$
happens to be symplectic if and only if $c_1=0$ and is
diffeomorphic $S^2 \times T^2$; furthermore we have that $L \cap
\mathbb R_+ \cdot F \cong \mathbb R \backslash \{0\}$ has a
homotopy type of a two--point set.

Given a closed nondegenerate 1--form on $N$ satisfying \ref{fro}
it is relatively easily to construct an invariant symplectic form
on $M$ such that $\iota_X\om=\pi^*\al.$ The construction follows
from a formula \cite{Bo} in the spirit of the {\it inflation}
trick of Thurston \cite{TH1} and McDuff \cite{MD}. It consists in
enlarging the form along the foliation determined by $\ker\alpha.$
For a given form $\alpha,$

\begin{eqnarray}\label{sim} \omega = \eta \wedge \pi^*\alpha + \pi^*(K \beta +
\phi)\end{eqnarray}

\noindent is an invariant symplectic form if $\eta \in
\Omega^1(M^4,\mathbb R)$ is a connection form, $\beta$ is a closed
2 -- form on $N^3$ such that $\alpha \wedge \beta$ is a volume
form on $N^3, \ d\phi = -c_1 \wedge \alpha$ and $K$ is
sufficiently large real number. Obviously $\omega$ satisfies
$\pi^*\alpha = \iota_X\omega$.

\smallskip

Existence of $\beta$ is well--known (\cite{Pl,Su}).

\begin{lemma}\label{bet} Let $N^n$ be a closed and oriented manifold. Assume
that a closed and nondegenerate 1--form $\alpha$ on $N$ is given.
Then there is a closed $(n-1)$--form $\beta$ such that $\alpha
\wedge \beta$ is a volume form on $M$. Equivalently, $\beta$ is
nondegenerate on leaves of the foliation defined by $\ker\alpha.$
\end{lemma}

\bigskip

\bigskip

We are able to state and prove the main theorem of this section.

\begin{thm}\label{gh} The total space $E$ of the fibration $T^2 \ra E {\buildrel \pi
\over \longrightarrow} \Sigma_g \ \ (g>1)$ is symplectic if and
only if $[T^2] \neq 0$ in $H_2(E,\mathbb R).$ \end{thm}

{\bf Proof.} Note that due to Theorem \ref{thc} it is enough to
prove that if $E$ is symplectic, then $[T^2] \neq 0.$ In the
sequel we shall need also the fact that $b_2(E)=2b_1(E)-2.$ This
is a direct consequence of the fact that the Euler characteristics
$\chi(E)$ vanishes. We will divide the proof into two parts:

\begin{enumerate}
\item bundles which does not admit a free circle action preserving
fibers, \item bundles which support such action.
\end{enumerate}

The first part will be proved by using spectral sequences. We
begin with flat $E$ bundle to obtain $b_2(E) = 2 +
rank(E_{11}^\infty)$. We change the Euler class keeping the
monodromy fixed to get $E'.$ Since $b_2(E)=b_2(E')$ and
$rank(E_{11}^\infty)=rank({E'}_{11}^\infty)$ we have $b_2(E') = 2
+ rank({E'}_{11}^\infty).$ But this means that all bundles in this
category have $\pi$--compatible symplectic forms.

The second part is proved in different way. Using section \ref{tp}
we know that if $E$ symplectic, then $E \slash S^1$ is the product
$\Sigma_g \times S^1.$ Furthermore, $E$ cannot be nontrivial
principal $T^2$--fibration over $\Sigma_g$ in order to be
symplectic (\cite{E}). However, it means that $L \neq A$
(\cite{HW}) and $E$ admits an invariant symplectic form. We also
get that for these bundles we have $[T^2] \neq 0.$

\subsection{First case} \label{first}

In the Leray--Serre spectral sequence for the fibration $T^2 \ra E
{\buildrel \pi \over \longrightarrow} \Sigma_g$ in the
$E^\infty$-- page we have

\begin{eqnarray}\label{ent} b_2(E)=rank(E_{20}^\infty) +
rank(E_{11}^\infty) + rank(E_{02}^\infty). \end{eqnarray}

\NI In addition to that, we know from Corollary \ref{fl} that
$[T^2]$ is nonzero in $H_2(E,\mathbb Z)$ since the bundle is flat.

\NI Finally lemma \ref{sps} simplifies the formula \ref{ent} to
$$ b_2(E)=rank(E_{11}^\infty) + 2.$$

\NI Note also the summand $E_{11}^\infty$ does not depend on the
Euler class, but only on monodromy.

Consider now the bundle $E'$ with the same monodromy as in the
bundle $E$ but nonzero Euler class. Lemma \ref{id} yields
$b_1(E')=2g=b_1(E)$ and $b_2(E')=2b_1(E')-2=2b_1(E)-2=b_2(E).$ In
both cases $rank(E_{11}^\infty)$ is the same, so
$b_2(E')=2+rank({E'}_{11}^\infty).$ But the latter means (Lemma
\ref{sps}) that $[T^2] \neq 0.$

\subsection{Second case}

We have here the total space $E$ equipped with a free circle
action preserving fibers. Furthermore, the quotient space $E
\slash S^1$ is a principal $S^1$--bundle over $\Sigma_g.$ From
section \ref{tp} we know that if $E$ is symplectic, then the
latter must be trivial and $E \slash S^1 \cong \Sigma_g \times
S^1.$ Therefore we shall restrict our considerations to bundles of
the form

\begin{eqnarray}\label{xii} \xi: S^1 \ra E \ra  \Sigma_g
\times S^1. \end{eqnarray}

Note that by \cite{E} we can also exclude principal nontrivial
$T^2$ fibrations over $\Sigma_g$ from our consideration.

We shall describe the Thurston norm on $\Sigma_g \times S^1.$ From
\cite{TH2} we know that ${\Arrowvert [\Sigma_g]\Arrowvert}_T
=|\chi(\Sigma_g)|=2g-2.$ Furthermore, if we denote the bundle
tangent to $\Sigma_g$ by $\tau$ and its Euler class by
$\chi(\tau)$ we have that the equality ${\Arrowvert a
\Arrowvert}_T =|\chi(\tau) \cdot a|$ holds for all $a \in
H_2(\Sigma_g \times S^1, \mathbb R)$ in some neighborhood of the
ray through $[\Sigma_g].$ In addition to that we know that second
homology classes of the form $b_i \times S^1, \ i=1,\ldots, 2g$
can be represented by embedded tori, where $b_i \hookrightarrow
\Sigma_g$ are embedded circles representing the base of
$H_1(\Sigma_g,\mathbb Z).$ It follows that the Thurston norm
vanishes identically on the subspace generated by classes of the
form $b_i \times S^1, \ i=1,\ldots, 2g.$ These data completely
describe the Thurston norm on $\Sigma_g \times S^1.$

When we pass to first real cohomology via Poincare duality we can
rewrite this data as ${\Arrowvert
k[S^1]+b\Arrowvert}_T=|k|(2g-2),$ where $b$ is any cohomology
class coming from $\Sigma_g.$ Thus all cohomology classes, except
of those which are pullbacks from $\Sigma_g,$ have a closed and
nondegenerate representative.

Denote then by $A \subset H^1(\Sigma_g \times S^1, \mathbb R)$ a
codimension 1 subspace generated by all these classes which are
pullbacks from $\Sigma_g.$ This subspace consists of classes on
which Thurston norm vanishes identically.

\smallskip

Recall that there exists an invariant symplectic form if and only
if $L \neq A$ (\cite{HW}). Thus $L$ (see \ref{na}) coincides with
$A$ if and only if $\xi$ is a nontrivial $T^2$--principal bundle
over $\Sigma_g.$ This bundle does not have any symplectic form by
\cite{E}. Recall that there exists an invariant symplectic form on
$E$ provided that $L \neq A$ (see \ref{sim}). In all these cases
$E$ supports an invariant symplectic form.

\smallskip

As we show below, whenever such $E$ supports an invariant
symplectic form, the homology class $[T^2]$ is nonzero in
$H_2(E,\mathbb R).$

\smallskip

Formula \ref{sim} gives the following application: if we have any
embedded circle $S^1$ in $N$ such that $\int_{S^1}\alpha>0,$ then
\begin{eqnarray}\label{com} \int_{\pi^{-1}(S^1)}\om = 2\pi
\int_{S^1}\alpha>0. \end{eqnarray}

Now, since $L$ does not coincide with $A$ we have that for any
nondegenerate 1--form $\al$ such that $[\al] \in L$ the inequality
$\int_{S'}\alpha>0$ holds, where $S'$ is the circle $S^1
\hookrightarrow S^1 \times \Sigma_g$ which is the first factor of
the product $S^1 \times \Sigma_g.$ But this in turn means that
$[T^2] \neq 0.$

\subsection{A characterization of the property $[T^2] \neq 0$}
\label{last}

In this section we give a complete characterization of the bundles
with nonzero homology class of the fiber in terms of monodromy and
the Euler class. We also give an alternative proof of the fact
that if $E$ has a free circle action preserving fibers, then the
relation $mx_1+nx_2$ belongs to the group generated by the
homology class of the orbit \ifff $[T^2] \neq 0$.

We start with bundles which do not admit a free circle action
preserving fibers (i.e. $\not \exists z \ A_iz=z,\  i=1,\ldots,
2g$). From Lemma (\ref{id}) and subsection \ref{first} we know
that for these bundles we always have $[T^2] \neq 0.$ By
Proposition \ref{per} the condition of admitting such action can
be described in terms of monodromy.

Consider now bundles which admit such free circle action. If such
bundle is flat, then $[T^2] \neq 0$ (see Corollary \ref{fl}). Any
other bundle $E'$ can be obtained by changing the Euler class and
keeping the monodromy fixed.

In the sequel we use notation of Lemma \ref{id}. From the homology
computations \ref{fi1}, \ref{se} we see that the first Betti
number can decrease by one. Recall that the subspace
$$S=\langle A_ix_1-x_1, A_ix_2-x_2,
i=1,2,\cdots,2g \rangle$$ has dimensional one over $\mathbb Z$
provided that $E$ is flat and nontrivial (Lemma \ref{id}). Take
also as before any basis $\{z,u\} \in \mathbb Z^2$ for which all
the monodromy matrices have the form
$\left(\begin{array}{cc} 1 & * \\
0 & 1 \end{array}\right).$ Obviously we have that $S<\mathbb Zz.$

If the extra relations $mx_1+nx_2 \in \mathbb Zz,$ then
$b_1(E)=b_1(E')$ and $b_2(E)=b_2(E').$ Following subsection
\ref{first} we get that $[T^2] \neq 0$ in this case.

Otherwise we have $mx_1+nx_2 \not\in \mathbb Zz.$ It yields
$b_1(E')=b_1(E)-1$ and $b_2(E')=b_2(E)-2.$ But this means that
$rank({E'}_{02}^\infty)=rank({E'}_{20}^\infty)=0.$ However the
last condition implies that $[T^2] = 0$ (see Lemma \ref{sps}).

\bigskip

These computations are concordant with our previous results. To
explain this consider a flat bundle $T^2 \ra E {\buildrel \pi
\over \longrightarrow} \Sigma_g$ with monodromy
of the form $\left(\begin{array}{cc} 1 & k_i \\
0 & 1 \end{array}\right)$ in some basis $\{z,u\}.$

This bundle is a $S^1$ bundle $\xi$ over the product $S^1 \times
\Sigma_g.$ This can be proven as follows. Each flat $T^2$--bundle
over $\Sigma_g$ has an obvious section, which assigns the neutral
element group $1 \in T^2$ to any point in base. This section
defines another section $\sigma: \Sigma_g \ra E\slash S^1,$ but
since $E\slash S^1$ is a principal $S^1$--bundle over $\Sigma_g,$
it must be trivial. As in \ref{xii} denote this bundle by $\xi$:

$$\xi: S^1 \ra E \ra  \Sigma_g \times S^1.$$

Using the Poincare duality it is easy to describe the first Chern
class of this bundle. Namely, it is equal to union of disjoint
circles representing the standard base $b_1, \ldots ,b_{2g} \in
H_1(\Sigma_g, \mathbb Z)$ with multiplicities equal accordingly to
$k_i.$

We next change the Euler class to obtain another bundle $E'.$ If
we have $mu+nv \in \mathbb Zu,$ then the section $\sigma$ is still
well defined as a mapping from $\Sigma_g$ to $T^2 \slash [z].$ The
latter means that the new bundle $E'$ can be interpreted as
another $S^1$ bundle $\xi'$ over the product $\Sigma_g \times
S^1:$

$$\xi': S^1 \ra E \ra  \Sigma_g \times S^1.$$

The difference between $\xi$ and $\xi'$ can be described by use of
the first Chern class. To be more precise: we have
$c_1(\xi)-c_1(\xi') \in \mathbb Z[S^1],$ where this circle is the
second factor in $\Sigma_g \times S^1.$

As we mentioned earlier, the argument concerning spectral
sequences gives that if $mu+nv \in \mathbb Zu,$ then $[T^2] \neq
0.$ It follows (Theorem \ref{thc}) that $E'$ supports a
$\pi$--compatible symplectic form. This form may be chosen to be
invariant, as the next lemma shows.

\begin{lemma}\label{av} If $E$ is equipped with a free circle action
included in fibers and $[T^2] \neq 0$, then $E$ has an invariant
$\pi$--compatible symplectic form.\end{lemma}

{\bf Proof.} It follows (\cite{TH1}) that $E$ supports a
symplectic structure $\om$ compatible with the fibration. The
symplectic form $\om$ can written down as

\begin{eqnarray}\om= \om_F+K\pi^*\om_B, \end{eqnarray}

\NI where $\om_B$ is the symplectic structure on the base
$\Sigma_g,$ $\om_F$ is a closed 2--form nondegenerate on the
fibers and $K$ is a big enough real number. Observe that $\om_F$
can be chosen invariant since we can average this form with
respect to the action. Finally, we can choose $K$ so big that the
form $\om$ is symplectic.

\QED

\begin{rem}\label{avr} {\em For torus bundles satisfying conditions that
$[T^2] \neq 0$ and $E$ is equipped with a free circle action
preserving fibers Lemma \ref{av} gives a simpler argument that
$b_1(E)\geq 2g+1$ (cf. Lemma \ref{id}). To see this denote by
$i:T^2 \subset E$ the inclusion of the fiber in the total space.
Next, take any basis $\{z,u\}$ for $H_1(T^2,\mathbb Z),$ where
$i_*z$ denotes the homology class of the orbit, and the invariant
$\pi$--compatible symplectic form $\om.$ Thus we obtain
$$0 \neq \int_u i^*(\iota_Z \om )=\int_{i_*u}\iota_Z \om,$$ where $Z$
denotes the infinitesimal generator of the action. This equality
yields that $i_*u \neq 0$ in $H_1(E,\mathbb R).$

Using this approach we can also provide another proof of Corollary
\ref{pr} (i.e. for every nontrivial principal torus fibration over
$\Sigma_g$ we have $[T^2] =0).$ Assume the opposite and take any
$\pi$--compatible symplectic form $\om.$ Similarly as Lemma
\ref{av} we can assume that $\om$ is averaged with respect to
principal torus action. Furthermore we have $\int_{i_*u}\iota_Z\om
\neq 0$ and $\int_{i_*z}\iota_U\om \neq 0,$ where $\{z,u\}$ is
some basis for $H_1(T^2,\mathbb Z)$ and $Z,U$ are corresponding
infinitesimal generators. However this means that $b_1(E)=2g+2$
(see formula \ref{se}) which contradicts Lemma \ref{id}.} \hfill
$\Box{}$
\end{rem}

\begin{rem}\label{avrr} {\em Averaging of forms can also be used to give a
simple proof of the following fact. Let $\xi:T^2 \ra E \ra M$ be a
principal torus bundle over some manifold $M.$ Then $[T^2] \neq 0$
in $H_2(E,\mathbb R)$ \ifff $\xi$ trivial. To see this assume
$[T^2] \neq 0.$ A slight generalization of the Thurston theorem
\ref{thc} yields existence of some closed two form $\beta$ on $E$
which is the volume form on each fiber. Average the form $\beta$
with respect to the torus action to get a closed invariant 2--form
which is also the volume form on each fiber. Now the same argument
as in Remark \ref{avr} gives that $u,v \neq 0$ in $H_1(E,\mathbb
R),$ where $\{u,v\}$ is some basis for $H_1(F,\mathbb R).$ This
yields $b_1(E)=b_1(M)+2$ which means that the bundle $\xi$ is
trivial. } \hfill $\Box{}$
\end{rem}

\bigskip

In this way we are able to answer which $T^2$ bundles over
$\Sigma_g$ has the property that $[T^2] \neq 0.$ They are
precisely:

\begin{enumerate}
\item $E=T^2 \times \Sigma_g$, \item bundles which do not support
any free circle action included in fiber, \item bundles which
support some free circle action included in fiber and the Euler
class is the multiple of the orbit class $[S^1] \in H_1(E,\mathbb
Z).$
\end{enumerate}

\newcommand{\gsv}{\mbox{$\tilde{\GG}_s(V)$}}
\newcommand{\gsw}{\mbox{$\tilde{\GG}_s(W)$}}
\newcommand{\gso}{\mbox{$\tilde{\GG}_s(0)$}}
\newcommand{\gsvw}{\mbox{$\tilde{\GG}_s(V\times W)$}}

\bibliographystyle{amsalpha}

\medskip
\noindent {\bf Mathematical Institute, Wroc\l aw University,

\noindent pl. Grunwaldzki 2/4,

\noindent 50-384 Wroc\l aw, Poland}

\begin{flushleft}
\tt rwalc@hera.math.uni.wroc.pl
\end{flushleft}

\end{document}